\documentclass[twocolumn]{autarthal}    

\usepackage{amsfonts}
\usepackage{amsmath}
\usepackage{amssymb}
\usepackage{graphicx}%
\usepackage[latin1]{inputenc}
\usepackage{epstopdf}
\usepackage{color, enumerate}
\usepackage{multicol}
\usepackage{subcaption}
\newcommand{\ben}{\begin{enumerate}}
	\newcommand{\een}{\end{enumerate}}
\newcommand{\R}{\mathbb{R}}
\newcommand{\N}{\mathbb{N}}
\newcommand{\blem}{\begin{lem}}
	\newcommand{\elem}{\end{lem}}
\newcommand{\bdfn}{\begin{dfn}}
	\newcommand{\edfn}{\end{dfn}}
\newcommand{\bcor}{\begin{cor}}
	\newcommand{\ecor}{\end{cor}}
\newcommand{\bthm}{\begin{thm}}
	\newcommand{\ethm}{\end{thm}}
\newcommand{\bex}{\begin{expl}}
	\newcommand{\eex}{\end{expl}}
\newcommand{\brmq}{\begin{rmq}}
	\newcommand{\ermq}{\end{rmq}}
\newcommand{\bitem}{\begin{itemize}}
	\newcommand{\eitem}{\end{itemize}}

\usepackage{graphicx}          

\begin{document}
	
	\begin{frontmatter}
		\runtitle{Event-triggered controls in abstract settings}  
		
		\title{Event triggered control and exponential stability for infinite dimensional linear systems \thanksref{footnoteinfo}} 
		
		\thanks[footnoteinfo]{
			S. E. is partially supported by the ANR projects TRECOS ANR 20-CE40-0009,  NumOpTes ANR-22-CE46-0005, CHAT ANR-24-CE40-5470.}
				
		\author[a]{Lucie Baudouin}\ead{lbaudouin@laas.fr},               
		\author[b]{Sylvain Ervedoza}\ead{sylvain.ervedoza@math.u-bordeaux.fr}  
		
		\address[a]{LAAS-CNRS, Universit\'e de Toulouse, CNRS, Toulouse, France}  
		\address[b]{Institut de Mathématiques de Bordeaux, UMR 5251, Université de Bordeaux, CNRS,
Bordeaux INP, F-33400 Talence, France.}  



%
%
%
%
%
%
\begin{abstract}
	This article aims at providing a unified analysis of the exponential stabilization of some abstract infinite dimensional systems undergoing an event-triggering mechanism that samples the control input. The partial differential equation is supposed to be defined by a skew-adjoint operator and controlled and observed through bounded operators. The continuously controlled closed loop system is assumed to be exponentially stable and the goal is to prove that a well-designed event-triggering mechanism to rule the time updates of the sampled control will allow to keep such a stability property.
	The key of the proof relies on the existence of an adequate Lyapunov functional. Existence and regularity of the solution to the closed-loop event-triggered system are also proven, along with the avoidance of Zeno behavior. 
\end{abstract}
	\end{frontmatter}

%
%
%
%

\section{Introduction}
Event-triggered control is a ``sample and hold" strategy where the update of the control in a closed loop system is triggered only when specific conditions are met.
Such strategies have emerged as a promising paradigm for the efficient use of computational resources in the control of dynamical systems. At a more fundamental level, this paper aims at providing a unified analysis in the context of abstract linear infinite dimensional systems, encompassing several models of  partial differential equations (PDEs), undergoing an event-triggering mechanism that samples the control input. 

Analyzing systems with event-based sampling poses specific challenges compared to those with periodic sampling due to the inherent time-varying nature of the closed-loop system. Extensive research has been conducted in the context of finite-dimensional systems of ordinary differential equations (ODEs), addressing various questions in seminal works \cite{Aarzen-Ifac99}, \cite{Tabuada-TAC07}, \cite{PANT-CDC11} and in more recent contributions such as \cite{PTNA-TAC15} (for nonlinear systems), \cite{TANWANI201646} (incorporating observers), \cite{PANS-TAC18} (employing high gain approaches),  \cite{Girard-TAC15} (with dynamic triggering mechanisms) and related references. 

Concerning infinite-dimensional systems, the literature about event-based control strategies is more recent and raising, e.g. for first order hyperbolic equations \cite{espitia2016event,EspAuriolYuKrstic-IJRNC22}, parabolic \cite{KFS-CDC19,SELIVANOV2016344,KEK-SCL24} or wave equations \cite{KBT-Auto22}, and dispersive equations like Korteweg-de Vries \cite{KBF} or Schrödinger \cite{KBT-SchroECC22}. One can also find now event-triggered control for systems coupling PDEs and ODEs, such as in \cite{WangKrstic-TAC23} (parabolic PDE-ODE cascade) or \cite{HSMdT-Auto25} (linearised FitzHugh-Nagumo system). This last article also offers a interesting overview of the state-of-the-art in event-triggered control for PDEs, paying attention to distinguish the recent achievements in hyperbolic and parabolic PDEs, since they correspond to very different phenomenons and have different properties that raise specific difficulties. We should also mention the results presented for infinite dimensional abstract settings in \cite{WakaikiSano-SICON20} about event triggered control (see also  \cite{LogRebTown} for sampled data controllers) which require an assumption of compactness of the feedback operator, that will not be compatible with the abstract PDE setting that we will be studying here.

PDEs, even restricted to linear ones, constitute indeed a rich mathematical framework for modeling very diverse physical phenomena, including for instance structural mechanics or vibrations, transport, or simplified models of fluid dynamics. However, controlling such systems presents significant challenges due to their infinite-dimensional nature and complex dynamics. 

In this article, we take advantage of some recent developments in event-triggered control approaches tailored for the wave \cite{KBT-Auto22} and Schr\"odinger \cite{KBT-SchroECC22} equations, and we build an abstract framework which encompasses, more broadly than these two examples, several models of interest that can be multidimensional in the space variable.  

The article will explore the event-triggering control strategy for the case of an abstract PDE system in an Hilbert setting in which the operator is a skew-adjoint infinitesimal generator of a $C_0$-semigroup, the closed loop is defined in continuous time by bounded control and observation operators and is exponentially stable. The usual questions will be answered about the new corresponding event-triggered closed loop system:  well-posedness in the appropriate functional space (to ensure the existence and uniqueness of solutions);  avoidance of Zeno phenomenon (where an infinite number of updates could be triggered in a finite time window); exponential stability of the state towards the same equilibrium.

Let us now dive into the mathematical framework and present our result more precisely. We first consider a control system taking the shape
\begin{equation}\label{System}
\left\{  
\begin{array}{lr}
\dot z (t) = Az(t) + Bu(t), \quad&  \qquad t\geqslant 0,\\
y(t) = Cz(t),& t\geqslant 0, \\
u (t)= Ky(t),& t\geqslant 0, \\
z(0) = z_0, 
\end{array}
\right. 
\end{equation}
in which no event-triggering mechanism appears. Here, the operator $A: \mathcal D(A) \subset H \mapsto H$ is the infinitesimal generator of a $C_0$-semigroup on a Hilbert space $H$ where $ \mathcal D(A)$ is the domain of $A$, $B$ is a control operator, $C$ is an observation operator and  $K$ is the control gain. 

In this article, we will work with the following standing assumptions:

{\bf Assumption 1:} The operator $A$ is skew adjoint on $H$.
\smallskip

{\bf Assumption 2:} $B\in \mathcal L(U, H)$, $C \in \mathcal L(H, Y)$, and $K \in \mathcal{L}(Y,U)$ and $U$ and $Y$ are assumed to be Hilbert spaces.
\smallskip

{\bf Assumption 3:} The operator $A+ BKC$ is the infinitesimal generator of an exponentially stable $C_0$-semigroup on $H$, namely
$$
	\exists M>0, \exists \alpha >0, \forall t \geqslant 0, \ \|e^{t(A+BKC)}\|_{\mathcal L(H)}\leqslant M e^{-\alpha t}.
\smallskip$$

In a nutshell, according to these assumptions, the operator $A$ generates a unitary group, the control, observation and gain operators (respectively, $B$, $C$ and $K$) are bounded, and the original closed loop system \eqref{System} is exponentially stable. Note that this is not possible if one of the operators $B$, $C$ or $K$ is compact, and thus our setting is not covered by the abstract framework proposed in \cite{WakaikiSano-SICON20}.
%
%
%
%

We are interested in investigating the situation where a pre-existing continuous feedback is modified via some aperiodic sample-and-hold mechanism.
More precisely, our objective is the study of the stabilisation of this abstract system when it is subjected to time event-triggered updates of the control $u$, 
i.e. replacing the control law $u(t) = Ky(t)$ by the following piecewise constant control law
\begin{equation}\label{discretecontrol}
	u(t) = Ky(t_k)~~~ \hbox{for all  }t \in [t_k,t_{k+1}), 
\end{equation}
where the time updates $(t_k)_{k \in \N}$ form an increasing sequence of times suitably chosen. 

In particular, our goal is to propose a well-chosen rule defining somehow appropriate time updates $t_k \in \mathbb R_+$ of the control $u$, that maintains the exponential stability of the closed-loop system. 
More precisely, in the spirit of \cite{KBT-SchroECC22,KBT-Auto22}, let us define the event-triggered law:
\begin{multline} \label{ETL}
	t_{k+1} = \sup \Big\{ t>t_k, 
	\forall \tau \in [t_k,t), \\ 
	\left\| C\big(z(\tau) - z(t_k)\big) \right\|^2_{Y} \leqslant \gamma^2 \| z(\tau)\|_H^2\Big\}.
\end{multline}
Within this law, the occurrence of an update is triggered by a tunable ratio $\gamma^2$ of the deviation between the measurement of natural evolution of the controlled state $C z(t)$ and the one of the steady state corresponding to the previous sampling instant $C z(t_k) $, with respect to a proportion of the current energy of the system $\| z(t)\|^2_H$.

%

Our goal is thus to study the dynamics of the system 
\begin{equation}\label{System-Trig}
\left\{  
\begin{array}{lr}
\dot z (t) = Az(t) + Bu(t), \quad&  \qquad t\geqslant 0,\\
y(t) = Cz(t),& t\geqslant 0, \\
z(0) = z_0, 
\end{array}
\right. 
\end{equation}
where the control $u$ is the piecewise constant control law given by the formula \eqref{discretecontrol} on each time interval $[t_k, t_{k+1})$ and the time updates $(t_k)_{k \in \N}$ are given by the event-triggered law \eqref{ETL}. In short, we will refer to this system as \eqref{discretecontrol}--\eqref{System-Trig}.
%
%


\noindent Therefore, three main specific questions are related to our study: 
\begin{itemize}
\item Prove the well-posedness of the event-triggered closed loop system \eqref{discretecontrol}--\eqref{System-Trig};
\item Check that the event-triggered law \eqref{ETL} does not generate Zeno phenomenon, where some accumulation point in the sequence $(t_k)$ could appear in finite time; 
\item Prove the exponential stability of the event-triggered closed loop system \eqref{discretecontrol}--\eqref{System-Trig}.
\end{itemize}

The core idea of this paper is to lean on the existence of a Lyapunov functional for the exponentially stable system \eqref{System} (the one without event-triggering mechanisms), presented in Section \ref{Sec-Lyapunov}, to analyze the stability of the system \eqref{System-Trig} with the control law given by \eqref{discretecontrol}-\eqref{ETL}. One of the keypoint of our analysis is the proof of the fact that there is no Zeno phenomenon, which will stem from a specific lemma relying on the boundedness of the control, observation and gain operators, and on the skew-adjointness of $A$. We refer to Section \ref{Sec-Well-Posedness-and-Stab-with-Trig} for detailed statements and proofs.

\section{A Lyapunov functional}\label{Sec-Lyapunov}
In this section, we consider a generic infinite dimensional system driven by the generator $\mathcal A$ of a $C_0$-semigroup on a Hilbert space $\mathcal{H}$:
\begin{equation}\label{SystemA}
\left\{  
\begin{array}{lr}
\dot z (t) = \mathcal Az(t) , \quad& t\geqslant 0,\\
z(0) = z_0.
\end{array}
\right. 
\end{equation}

Our first result is somewhat classical, establishing the existence of a Lyapunov functional on $\mathcal{H}$ if the $C_0$-semigroup $(e^{t\mathcal A})_{t\geqslant 0}$ is exponentially stable:

\begin{thm}\label{LYAPU}
Let the operator $\mathcal A: \mathcal D(\mathcal A) \subset \mathcal{H} \mapsto \mathcal{H}$ be the infinitesimal generator of a  $C_0$-semigroup $\big(e^{t\mathcal A}\big)_{t\geqslant 0}$. 
Assume that there exists $c_1\geqslant 0$ such that 
\begin{equation}\label{hypmax}
	\forall z\in \mathcal D(\mathcal A),\quad \langle z, \mathcal A z \rangle_{\mathcal{H}} \leqslant c_1 \|z\|_{\mathcal{H}}^2.
\end{equation}
Assume also that the  $C_0$-semigroup $\big(e^{t\mathcal A}\big)_{t\geqslant 0}$ is exponentially stable, 
meaning that there exist $M>0$ and  $\alpha>0$ such that, for all $z_0 \in H$ and for all $t \geqslant 0$, 
\begin{equation}\label{hypexp}
 \left\|z(t) \right\|_{\mathcal{H}} = \left\|e^{t\mathcal A} z_0 \right\|_{\mathcal{H}} \leqslant M e^{-\alpha t} \|z_0\|_{\mathcal{H}}.
\end{equation}
Then for any $\beta \in (0, \alpha)$, there exists  $V_\beta : \mathcal{H} \to \mathbb R_+^*$ satisfying 
\begin{enumerate}
\item[(i)] $V_\beta \sim \|\cdot\|_{\mathcal{H}}^2$, meaning that ~ $\exists C_1,C_2>0$, $\forall z_0 \in \mathcal{H}$, 
$$C_1  \|z_0\|_{\mathcal{H}}^2 \leqslant V_\beta(z_0) \leqslant  C_2  \|z_0\|_{\mathcal{H}}^2 ~;$$

\item[(ii)] For any $z_0\in \mathcal{D}( \mathcal A)$, the solution $z$ of \eqref{SystemA} satisfies, for all $t \geqslant 0$, 
$$\dfrac d{dt} V_\beta(z(t)) \leqslant -2\beta V_\beta(z(t)), $$ 
 or, equivalently, for any $z_0\in \mathcal{D}( \mathcal A)$,
 $$ \langle \nabla V_\beta(z_0) , \mathcal A z_0\rangle_{\mathcal{H}} \leqslant - 2 \beta V_\beta(z_0) ~; $$ 

\item[(iii)] There exists $C_0 > 0$ such that for all $ z_0 \in \mathcal{H}$, 
$$\| \nabla V_\beta(z_0)\| \leqslant C_0 V_\beta(z_0)^{\frac 12}.$$
\end{enumerate}
Choosing the Lyapunov functional  $V_\beta$ suitably, the above constants can be chosen as 
\begin{multline}
	\label{Constants-eta}
	 C_0 = 2\sqrt{1+(c_1+\beta) M^2}, \quad \text{ and } 
			\\
	 C_1 = \frac{1}{2}, 
	\quad
	C_2 = \frac{1}{2} \left( 1+ \frac{(c_1+\beta) M^2}{\alpha - \beta}\right).
\end{multline}
\end{thm}


One can read a similar result in the books \cite[Theorem 5.1.3]{CurtainZwart} or  \cite[Theorem 8.1.3]{JZbook}, giving the equivalence between the exponential stability of a  $C_0$-semigroup and the existence of a positive operator, solution of the associate Lyapunov equation. Nevertheless, we decided to present here the statement and proof of the existence of a Lyapunov functional for an exponentially stable system in order to offer a self contained reasoning, and get an explicit grasp on the values of the constants $C_0$, $C_1$ and $C_2$ appearing in Theorem \ref{LYAPU}.\\

\begin{pf}
For $\beta \in (0, \alpha)$ and $\eta \geqslant 0$, we introduce the functional $V_{\beta, \eta}$ defined for $z_0\in \mathcal{H}$ by
\begin{equation}\label{Lyap}
V_{\beta, \eta} (z_0) = \dfrac 12 \|z_0\|_{\mathcal{H}}^2 + \eta \int_0^{+\infty} e^{2\beta s} \| e^{s\mathcal A} z_0\|_{\mathcal{H}}^2 \,ds.
\end{equation}
The first remark is that, thanks to the assumption \eqref{hypexp} and since $\beta < \alpha$, this functional is always well defined for any $z_0\in \mathcal{H}$. Next, we will check that for a suitable choice of $\eta$, the functional $V_{\beta, \eta}$ is a Lyapunov functional for system \eqref{SystemA}.

\noindent $\bullet$  First, it is obvious that for any $z_0\in \mathcal{H}$, one has $\frac 12 \|z_0\|_{\mathcal{H}}^2 \leqslant V_{\beta, \eta}(z_0) $. Moreover, the assumption \eqref{hypexp} and the fact that $\beta < \alpha$ bring that for all $z_0\in \mathcal{H}$, 
\begin{eqnarray*}
V_{\beta, \eta} (z_0)  &\leqslant&  \dfrac 12 \|z_0\|_{\mathcal{H}}^2 + \eta M^2 \int_0^{+\infty} e^{2\beta s} e^{-2\alpha s} \|z_0\|_{\mathcal{H}}^2 \,ds\\
&\leqslant&  \dfrac 12\left(1 + \dfrac{\eta M^2}{\alpha-\beta} \right)\|z_0\|_{\mathcal{H}}^2.
\end{eqnarray*}
Both sides of the estimate $(i)$ stating the equivalence of the proposed functional $V_{\beta, \eta}$ with the usual norm of the Hilbert space $\mathcal{H}$ are thus proved, whatever the value of $\eta \geq 0$.

\noindent $\bullet$  Second, one  calculates the derivative of $t\mapsto V_{\beta, \eta}(z(t))$ for any $z_0 \in \mathcal D(\mathcal A)$. Using \eqref{SystemA}, $z(t) = e^{t\mathcal A} z_0$, and an integration by parts, one gets, for $t \geqslant 0$, 
\begin{eqnarray*}
&\frac d{dt}& V_{\beta, \eta}(z(t))  \\
&=&  \langle z(t), \mathcal Az(t) \rangle_{\mathcal{H}} + \eta \int_0^{+\infty} e^{2\beta s} \frac \partial{\partial t}\left(  \| e^{(s+t)\mathcal A} z_0\|_{\mathcal{H}}^2\right)  \,ds \\
&=& \langle z(t), \mathcal Az(t) \rangle_{\mathcal{H}} + \eta \int_0^{+\infty}  e^{2\beta s} \frac \partial{\partial s}\left(  \| e^{(s+t)\mathcal A} z_0\|_{\mathcal{H}}^2\right)  \,ds \\
&=&  \langle z(t), \mathcal Az(t) \rangle_{\mathcal{H}} -2 \eta\beta \int_0^{+\infty}  e^{2\beta s}   \|e^{(s+t)\mathcal A} z_0\|_{\mathcal{H}}^2  \,ds \\
&&+~  \eta \left[ e^{2\beta s} \| e^{(s+t)\mathcal A} z_0\|_{\mathcal{H}}^2\right]_{s=0}^{s\to+\infty}.
\end{eqnarray*}
On the one hand, at $s=0$, we get 
$$\left.e^{2\beta s} \| e^{(s+t)\mathcal A} z_0\|_{\mathcal{H}}^2\right|_{s=0} \\
=  \| e^{t\mathcal A} z_0\|_{\mathcal{H}}^2 = \| z(t)\|_{\mathcal{H}}^2,
$$
and on the other hand, since $\beta < \alpha$, the exponential stability of the semigroup $\big(e^{s \mathcal A}\big)_{s \geqslant 0}$  brings 
\begin{multline*}
\lim_{s\to +\infty} \left( e^{2\beta s} \| e^{(s+t)\mathcal A} z_0\|_{\mathcal{H}}^2\right) \\
= \lim_{s\to +\infty} \left( e^{2\beta s} \| e^{s \mathcal A} (e^{t \mathcal A} z_0)\|_{\mathcal{H}}^2\right) \\
\leqslant  \lim_{s\to +\infty} \left( M^2 e^{2(\beta-\alpha) s} \|(e^{t \mathcal A} z_0)\|_{\mathcal{H}}^2\right) = 0.
\end{multline*}
Therefore, using \eqref{hypmax} and \eqref{Lyap}, one obtains, if $\eta  \geqslant c_1+\beta$,
\begin{eqnarray*}
	&\frac d{dt}&\left( V_{\beta, \eta}(z(t)) \right)
	 =   \langle z(t), \mathcal Az(t) \rangle_{\mathcal{H}} \\
	 &&\quad \quad -2 \eta\beta \int_0^{+\infty}  \hspace{-0.2cm} e^{2\beta s}   \|e^{(s+t)\mathcal A} z_0\|_{\mathcal{H}}^2  \,ds 
	 -  \eta \| z(t)\|_{\mathcal{H}}^2\\
	&\leqslant&  c_1 \|z(t)\|_{\mathcal{H}}^2  - 2 \beta V_{\beta, \eta}(z(t)) +\beta \|z(t)\|_{\mathcal{H}}^2 -  \eta \| z(t)\|_{\mathcal{H}}^2\\
	&\leqslant&  -2 \beta V_{\beta, \eta}(z(t)).
\end{eqnarray*}

The other statement of $(ii)$ can be proved on its own but also simply stems from the fact that for any $z_0 \in \mathcal D(\mathcal A)$, considering the solution $z$ of \eqref{SystemA}, we have 
\begin{eqnarray*}
	\left. \left(\frac d{dt}\left( V_{\beta, \eta} (z(t)) \right) \right)\right|_{t = 0}&=&  \left(\langle \nabla V_{\beta, \eta}( z(t)),  \dot z(t) \rangle_{\mathcal{H}}\right)|_{t = 0}
	\\ &=& \langle \nabla V_{\beta, \eta}( z_0), \mathcal Az_0 \rangle_{\mathcal{H}}.
\end{eqnarray*}

\noindent $\bullet$  Third, in order to prove $(iii)$, let us first calculate that, for any $ z_0 \in \mathcal{H}$, 
$$
	 \nabla V_{\beta, \eta}(z_0) = z_0 + 2\eta \int_0^{+\infty} e^{2\beta s} e^{s\mathcal A^*}e^{s\mathcal A}  z_0  \,ds.
$$ 
Then, we have, using $\beta < \alpha$ and the exponential stability of $\big(e^{t \mathcal A}\big)_{t \geqslant 0}$ 
(thus of  $\big(e^{t \mathcal A^*}\big)_{t \geqslant 0}$ with the same constants):
\begin{eqnarray*}
	\lefteqn{
	\| \nabla V_{\beta, \eta}(z_0)\|_{\mathcal{H}}^2 =  \langle \nabla V_{\beta,\eta}(z_0), \nabla V_{\beta, \eta}(z_0)\rangle_{\mathcal{H}}
	}\\
		\qquad &=& \| z_0\|_{\mathcal{H}}^2 + 4\eta    \int_0^{+\infty} e^{2\beta s} \|e^{s\mathcal A} z_0 \|_{\mathcal{H}}^2 \,ds \\
		&&+ ~4 \eta^2 \left\| \int_0^{+\infty} e^{2\beta s} e^{s\mathcal A^*} e^{s\mathcal A} z_0  \,ds \right\|_{\mathcal{H}}^2\\
		&\leqslant& \| z_0\|_{\mathcal{H}}^2 + 4\eta    \int_0^{+\infty} e^{2\beta s} \|e^{s\mathcal A} z_0 \|_{\mathcal{H}}^2 \,ds\\
		&&+ ~ 4 \eta^2  M^2 \int_0^{+\infty} e^{4\beta s} e^{-2\alpha s }\|e^{s\mathcal A} z_0 \|_{\mathcal{H}}^2 \,ds  \\
		&\leqslant& \| z_0\|^2 + 4\eta (1+\eta  M^2)   \int_0^{+\infty} e^{2\beta s} \|e^{s\mathcal A} z_0 \|_{\mathcal{H}}^2 \,ds \\
		&\leqslant& \max\left\{2\, ,4 (1+\eta  M^2)\right\} V_{\beta, \eta}(z_0)~ \leqslant ~C_0^2  V_{\beta, \eta}(z_0),
\end{eqnarray*}
with $ C_0 = 2\sqrt{1+\eta M^2}$. 

Gathering all the estimates proved above, we proved that for any $\beta \in (0,\alpha)$ and $\eta \geq c_1 + \beta$, $V_{\beta, \eta}$ is a suitable Lyapunov functional for \eqref{SystemA} satisfying all the properties in Theorem \ref{LYAPU}. Note that all the constants here corresponding to $V_{\beta, \eta}$ can be estimated explicitly as 
\begin{multline*}
	C_{0}(\eta) = 2\sqrt{1+\eta M^2}, 
	\quad \text{ and } 
	\\
	C_1(\eta) = \frac{1}{2}, 
	\quad 
	C_2(\eta) = \frac{1}{2} \left( 1+ \frac{\eta M^2}{\alpha - \beta}\right).
\end{multline*}
To have them as small as possible and obtain the closest possible Lyapunov functional to $\| \cdot \|_{\mathcal{H}}^2$, it is clear that one should take $\eta = c_1 + \beta$. We thus  set $V_\beta = V_{\beta, c_1 + \beta}$, for which the constants in Theorem \ref{LYAPU} can then be chosen as in \eqref{Constants-eta}.
$\hfill$ \qed
\end{pf}
\section{Event triggered stabilization}\label{Sec-Well-Posedness-and-Stab-with-Trig}
Let us now consider the closed loop system \eqref{discretecontrol}--\eqref{System-Trig}, that we rewrite in a more concise form 
\begin{equation}\label{SystemEch}
\left\{  
\begin{array}{ll}
\dot z (t) = Az(t) + Bu(t), \quad & t\in \mathbb R_+,\\
u (t)= KCz(t_k),\quad& t\in [t_k,t_{k+1}), k\in \mathbb N,\\
z(0) = z_0
\end{array}
\right. 
\end{equation}
submitted to the event triggering mechanism starting from $t_0 = 0$ and following the update law
\begin{multline} \label{ETL1}
	t_{k+1} = \sup \Big\{ t>t_k, 
	\forall \tau \in [t_k,t), \\ 
	\left\| Cz(\tau) - Cz(t_k) \right\|^2_{Y} \leqslant \gamma^2 \| z(\tau)\|_H^2\Big\}.
\end{multline}
where  $\gamma$ is a positive constant to be defined. 

Since the sampling is aperiodic, we will have to carefully check whether it does, or not, produce any Zeno phenomenon. And in that perspective, let us introduce
 the maximal time $T^*$ under which the closed loop system subjected to this event-triggered law has a solution:
\begin{equation}
\label{T}
\left\{
\begin{array}{ll}
T^*=+\infty  &  \text{ if $(t_k)$ is a finite sequence, }\\
T^*=\displaystyle\limsup_{k\to +\infty} t_k &  \text{ if not}.
\end{array}
\right.
\end{equation}
Hence, later in the article, the proof of the absence of Zeno behavior will actually result from the proof that $T^*=+\infty $, since no accumulation point of the sequence $(t_k)_{k\geqslant 0}$ will be possible. 

Note that the sequence of times $(t_k)_{k \in \N}$ and $T^*$ defined by \eqref{T} depend on the initial data $z_0$ through the dynamics \eqref{SystemEch}--\eqref{ETL1}. This dependence is not explicitly written to simplify the notations. 

\begin{thm}[Well-posedness up to $T^*$]\label{WellPo}
Under Assumptions 1 and 2, the closed loop system \eqref{SystemEch} submitted to the event-triggered law \eqref{ETL1} is well-posed in $H$ until the possibly infinite time $T^*$. More precisely, for $z_0 \in H$, there exist a unique triple $(z, (t_k)_{k\in \N}, T^*)$ with $z \in \mathcal{C}([0,T_*); H)$ and $(t_k)_{k \in \N}$ satisfying \eqref{ETL1} and \eqref{T} such that for all $t \in (0,T^*)$, 
$$
	z(t) = e^{t A} z_0 + \sum_{t_{k} < t} \int_{t_k}^{\min\{t, t_{k+1}\} } e^{ (t-s) A} BKC z(t_k) \, ds.
$$
Moreover, there exists a contant $\kappa>0$ such that for any $z_0 \in H$, the solution $z$ satisfies, for any $t\in[0,T^*)$,
\begin{equation}\label{estimWP}
 \|z_0\|_H e^{-\kappa t} \leqslant \|z(t)\|_H \leqslant  \|z_0\|_H e^{\kappa t}.
\end{equation}
Actually, one can take $\kappa = \|BKC\|_{\mathcal L(H)} + \gamma \|B K\|_{\mathcal L(Y,H)}$.
\end{thm}
\begin{pf} We will argue by induction and assume that the solution of \eqref{SystemEch}--\eqref{ETL1} is constructed up to the time $t_k$ for some $k \in \N$. 

We first have to check that $t_{k+1}$ given by \eqref{ETL1} is strictly larger than $t_k$. To do so, let us consider the system 
$$
\dot z_k(t) = Az_k(t) + BK C z(t_{k}), \text{ for } t \geqslant t_k, \quad z_k(t_k) = z(t_k).
$$
By Assumption 2, the operator $BKC$ belongs to $\mathcal L(H)$. Therefore, $BKC z(t_k) \in L^1_{loc}(t_k, \infty; H)$ and the solution $z_k$ is well-defined, for  $t \geqslant t_k$, by the Duhamel formula 
$$
	z_k(t) = e^{(t-t_k) A} z_k(t_k) + \int_{t_k}^t e^{(t -s ) A} BKC z(t_k) \, ds, 
$$
and belongs to $\mathcal{C}([t_k,\infty;H)$. \\
Note that, since $C$ belongs to $\mathcal L(H,Y)$, the deviation term $t \mapsto Cz_k(t) - Cz(t_k)$ is a continuous function on $[t_k,\infty)$, vanishing at $t=t_k$, so that $t_{k+1}$ is necessarily strictly larger than $t_k$, well-defined and unique.

Of course, the solution $z$ of \eqref{SystemEch}--\eqref{ETL1} is then simply given as the restriction of $z_k$ on the time interval $[t_k, t_{k+1})$.


To get estimate \eqref{estimWP}, note that we can also write the equation \eqref{SystemEch} on $(t_k, t_{k+1})$ under the form 
$$
\dot z (t) = (A + BKC) z(t) + B K C( z(t_k) - z(t)),
$$
so that one can calculate, using Assumption 1,
\begin{eqnarray*}
&\dfrac d{dt}&\left( \dfrac 12 \|z(t)\|_H^2\right) = \langle \dot z(t), z(t) \rangle_H\\
&=&\langle (A + BKC) z(t), z(t) \rangle_H + \langle B K C( z(t_k) - z(t)), z(t) \rangle_H \\
&=&\langle BKC z(t), z(t) \rangle_H + \langle B K C( z(t_k) - z(t)), z(t) \rangle_H.
\end{eqnarray*}
Using Assumption 2 and the event-triggered law \eqref{ETL1} that imply 
$\left\| C(z(t) - z(t_k)) \right\|_{Y} \leqslant \gamma \| z(t)\|_H$ for all $t\in [t_k,t_{k+1}[$, we get
\begin{multline*}
\left| \dfrac d{dt}\left( \|z(t)\|_H^2\right) \right| 
\leqslant ~2\|BKC\|_{\mathcal L(H)} \|z(t)\|_H^2 \\
+ 2\|B K\|_{\mathcal L(Y,H)} \|C( z(t_k) - z(t))\|_Y \|z(t) \|_H \\
\leqslant2\left(\|BKC\|_{\mathcal L(H)} + \gamma\|B K\|_{\mathcal L(Y,H)} \right) \|z(t) \|_H^2.
\end{multline*}
Therefore,  denoting $\kappa = \|BKC\|_{\mathcal L(H)} + \gamma \|B K\|_{\mathcal L(Y,H)}$, we get both
$$\dfrac d{dt}\left( e^{-2\kappa t}\|z(t)\|_H^2\right) \leqslant 0 \quad  \hbox{ and } \quad \dfrac d{dt}\left( e^{2\kappa t}\|z(t)\|_H^2\right) \geqslant 0$$ 
so that integrating over $(t_k,t)$, we obtain, $\forall t\in(t_k,t_{k+1})$,
$$\|z(t_k)\|_H^2 e^{-2\kappa(t-t_k)} \leqslant \|z(t)\|_H^2 \leqslant \|z(t_k)\|_H^2 e^{2\kappa(t-t_k)}.$$
Using a recurrence argument and the definition \eqref{T} of $T^*$, estimate \eqref{estimWP} and Theorem~\ref{WellPo}  are then easily obtained. 
$\hfill$ \qed 
\end{pf}

An important result to exclude the Zeno phenomenon is the following one:
\begin{prop}[More regularity]\label{MoReg}
Let us make  Assumptions~1 and 2. If $z_0\in \mathcal D(A)$, then the solution $z$ of system \eqref{SystemEch}--\eqref{ETL1} is weakly continuous in $[0,T^*)$ with values in $\mathcal{D}(A)$, its time derivative $\dot z$ is piecewise continuous in $[0,T^*)$ with values in $H$, with jumps at the discrete times $(t_k)_{k \in \N}$, and for all $t \in [0, T^*)$,
\begin{equation}\label{more}
 \|\dot z\|_{L^\infty(0,t;H)} \leqslant \|(A+BKC)z_0\|_H e^{t \|BKC\|_{\mathcal L(H)}}.
\end{equation}
\end{prop}
\begin{pf} Let $k\in \mathbb N$, and assume that $z(t_k) \in \mathcal D(A)$. 
Then $z$ solves
\begin{equation}\label{eqk}
\dot z (t) = Az(t) + BKCz(t_k),\quad t\in [t_k,t_{k+1}).
\end{equation}
It is then clear that the solution $ \zeta$ of 
$$
\left\{  
\begin{array}{ll}
\dot \zeta (t) = A\zeta(t), \quad& t\in [t_k,t_{k+1}),\\
\zeta(t_k) =  \left(A+ BKC\right)z(t_k)
\end{array}
\right. 
$$
satisfies $\zeta \in \mathcal{C}([t_k, t_{k+1}); H)$ since  $z(t_k) \in \mathcal D(A)$ and $BKC \in \mathcal L(H)$ by Assumption 2, and 
$$z(t) = z(t_k) + \int_{t_k}^t \zeta(s) \, ds.$$ 
Since $A$ is the generator of a unitary group on $H$ by Assumption 1, this yields that $\dot z \in \mathcal{C}([t_k, t_{k+1}); H)$, and for all $t \in [t_k, t_{k+1})$,
\begin{multline}\label{estdot}
	\|\dot z(t) \|_H = \|\zeta(t) \|_H 
	\\
	= \|\zeta(t_k) \|_H = \|  \left(A+ BKC\right)z(t_k) \|_H.
\end{multline}
Using \eqref{eqk} and Theorem \ref{WellPo}, we derive $Az \in \mathcal{C}([t_k, t_{k+1}); H)$, i.e. $z \in \mathcal{C}([t_k, t_{k+1}); \mathcal{D}(A))$, and that $\forall t \in [t_k,t_{k+1})$,
$$
	\|Az(t) + BKCz(t_k) \|_H  = \|  \left(A+ BKC\right)z(t_k) \|_H, 
$$
and 
\begin{equation}\label{bound}
	\| A z(t) \|_H \leq \|  \left(A+ BKC\right)z(t_k) \|_H + \| BKCz(t_k) \|_H.
\end{equation}
Accordingly, for $t \in [t_k, t_{k+1})$,
\begin{eqnarray*}
&\|&  \left(A+ BKC\right)z(t) \|_H
 \\
&\leqslant &\| Az(t) + BKC z(t_k)  \|_H + \| BKC\left(z(t) - z(t_k)\right) \|_H\\
&\leqslant & \| \left(A+ BKC\right)z(t_k) \|_H \\
&&+ \| BKC \|_{\mathcal L(H)}\|\dot z\|_{L^\infty(t_k,t;H)}(t- t_k)
\end{eqnarray*}
and obtain finally, using  \eqref{estdot}, that for $t \in [t_k, t_{k+1})$,
\begin{align}
\|&  \left(A+ BKC\right)z(t) \|_H \nonumber\\
& \leqslant  \| \left(A+ BKC\right)z(t_k) \|_H \left(1+ (t- t_k) \| BKC \|_{\mathcal L(H)} \right)\nonumber\\
&\leqslant \| \left(A+ BKC\right)z(t_k) \|_H ~e^{ (t- t_k)\| BKC \|_{\mathcal L(H)}}.\label{exp}
\end{align}
Since $z$ is continuous in $t_{k+1}$ with values in $H$ from Theorem \ref{WellPo}, passing to the limit $t \to t_{k+1}^-$ in estimates \eqref{bound} and \eqref{exp} brings that $Az(t)$ converges to $Az(t_{k+1})$ weakly in $H$ (hence the weak continuity of $z$ in time with values in $\mathcal{D}(A)$, but also $z(t_{k+1}) \in \mathcal{D}(A)$) and 
\begin{align*}
\|&  \left(A+ BKC\right)z(t_{k+1}) \|_H \\
&\leqslant \| \left(A+ BKC\right)z(t_k) \|_H ~e^{ (t_{k+1}- t_k)\| BKC \|_{\mathcal L(H)}}.
\end{align*}

By recurrence, we therefore prove that for any $k\in \mathbb N$,
$$
\|  \left(A+ BKC\right)z(t_k) \|_H \leqslant \| \left(A+ BKC\right)z_0 \|_H ~e^{ t_k\| BKC \|_{\mathcal L(H)}}.
$$
allowing to deduce from \eqref{estdot}  and the definition of $T^*$ that for all $t< T^*$
$$
 \|\dot z\|_{L^\infty(0,t;H)} \leqslant \|(A+BKC)z_0\|_H ~e^{t \|BKC\|_{\mathcal L(H)} }.
$$
The regularities on $z$ and $\dot z$ announced in Proposition \ref{MoReg} easily follow. $\hfill$ \qed 
\end{pf}

We are then in position to guarantee the absence of Zeno phenonemon:
\begin{thm}[No Zeno]\label{NoZe}
Under Assumptions 1 and 2,
if $z_0 \in \mathcal D(A)$, the solution $\big(z, (t_k)_{k\in \N}, T^*\big)$ of \eqref{SystemEch}-\eqref{ETL1} satisfies $T^* = +\infty$.\\
As a consequence, as soon as the initial state $z_0$ belongs to $ \mathcal D(A)$, the closed loop system \eqref{SystemEch}--\eqref{ETL1} do not experience any Zeno phenomenon.
\end{thm}
\begin{pf}	We argue by contradiction and assume that $T^*< + \infty$.
Since the sequence $(t_k)_{k\in \mathbb N}$ is increasing, from \eqref{T} we have $t_k \underset{k\to+\infty}{\longrightarrow} T^*$.\\
In the perspective of finding a contradiction, we aim at estimating $t_{k+1} - t_k$ by below for all $k \in \N$, as if we were looking for a ``dwell time" between two updates. We thus fix $k \in \N$, and for any $t\in (t_k,t_{k+1})$, we define 
$$\psi_k(t) = \|C(z(t)- z(t_k))\|^2_Y.$$ 
Let us recall that $\psi_k(t_k) = 0$, and that the event triggered law \eqref{ETL1} implies that $\psi_k(t) \leqslant \gamma^2 \|z(t)\|^2_H$ for all $t\in (t_k,t_{k+1})$, and $\lim_{t\to t_{k+1}}\psi_k(t) =  \gamma^2 \|z(t_{k+1})\|_H^2$.

The estimate we are looking for will stem from the study of the time derivative of this function $\psi_k$ on $(t_k, t_{k+1})$. Recalling that $z \in \mathcal{C}^1([t_k,t_{k+1}); H)$ for $z_0 \in \mathcal{D}(A)$, (cf. Proposition \ref{MoReg}), we get 
\begin{multline*}
\dot \psi_k(t)= 2 \langle C(z(t)- z(t_k)) , C \dot z(t) \rangle_Y \\
\leqslant  2 \gamma  \|z(t)\|_H \|C\|_{\mathcal L(H,Y)}\| \dot z(t) \|_H.
\end{multline*} 
Integrating $\psi_k$ over $(t_k,t_{k+1})$ brings
\begin{multline*}
 \gamma^2 \|z(t_{k+1})\|_H^2 \leqslant  2 \gamma  \|C\|_{\mathcal L(H,Y)}\int_{t_k}^{t_{k+1}} \|z(t)\|_H\| \dot z(t) \|_H\,dt\\
 \leqslant  2 \gamma  \|C\|_{\mathcal L(H,Y)} \left(t_{k+1}-t_k\right) \|\dot z\|_{L^\infty([0,T^*);H)}  \| z\|_{L^\infty([0,T^*);H)}.
\end{multline*} 
Therefore, using the lower estimate of \eqref{estimWP} from Theorem~\ref{WellPo} and  estimate \eqref{more} from Proposition~\ref{MoReg} we obtain
$$
t_{k+1} - t_k \geqslant \dfrac{\gamma e^{-(3\kappa + \|BKC\|_{\mathcal L(H)}) T^*}}{2   \|C\|_{\mathcal L(H,Y)}}  \dfrac{ \|z_0\|_H}{ \|(A+BKC)z_0\|_H  }.
$$
Since this estimate should be true for any $k \in \N$, this would bring a contradiction with $T^*< + \infty$.
$\hfill$ \qed 
\end{pf}

\begin{rem}
In the context of finite dimensional systems undergoing event-triggered control or in the study of hybrid dynamical systems, the proof that no Zeno phenomenon will occur often relies on the proof of the existence of a dwell-time $\tau$ between updates instants or jumps.  It is usually expected that this dwell-time $\tau$ doesn't depend on the time window $[0,T]$ on which the study is done, which is consistent with an engineering point of vue, and brings the avoidance of Zeno behavior. Nevertheless, recall that Zeno behavior is defined as an infinite number of jumps in finite time and therefore, it is equivalent with an accumulation point for the updates sequence. Here, we ruled out Zeno phenomenon proving that no such thing can occur. But our proof does not yield the existence of a positive dwell-time on the whole dynamics: it only yields that for all $z_0 \in \mathcal{D}(A)$, for all $T>0$, 
\begin{multline*}
	\inf_{t_{k+1} \leqslant T} \{ t_{k+1} - t_k \} 
	\\
	\geqslant \dfrac{\gamma e^{-(3\kappa + \|BKC\|_{\mathcal L(H)}) T}}{2   \|C\|_{\mathcal L(H,Y)}}  \dfrac{ \|z_0\|_H}{ \|(A+BKC)z_0\|_H  }, 
\end{multline*}
which goes to $0$ as $T \to \infty$.
\end{rem}
Let us now state the exponential stability result we can prove for the event-triggered closed loop system.
\begin{thm}[Event triggered stability]\label{ETES}
Under Assumptions 1, 2 and 3, the closed loop system \eqref{SystemEch} submitted to the event-triggered law \eqref{ETL1} is exponentially stable in $H$ for initial data in $\mathcal{D}(A)$ as soon as the tuning parameter $\gamma$ of \eqref{ETL1} is chosen small enough.  More precisely, setting the radius
$$r = \|BK\|_{\mathcal L(Y,H)}\sqrt{2(1+( \|BKC\|_{\mathcal L(H)}  + \alpha)  M^2)},$$ 
where $\alpha$ is the decay rate of the exponentially stable system \eqref{System} and $M$ its overshoot constant (in other words, $\alpha$ and $M$ are the constants appearing in Assumption 3),
then if we choose the tuning parameter $\gamma$ such that
\begin{equation}
	\label{Cond-SmallGamma} 
	\gamma < \dfrac{\alpha}{r }
\end{equation}
then for all $\delta >0$ satisfying 
\begin{equation}
	\label{Cond-delta}
	\delta < \alpha - \gamma r,
\end{equation}
there exists $C_* >0$ such that for $z_0 \in \mathcal{D}(A)$ the solution of \eqref{SystemEch} - \eqref{ETL1} satisfies
$$
\|z(t)\|_H \leqslant C_* \|z_0\|_H e^{- \delta t}, \qquad \forall t \geqslant 0.
$$
\end{thm}
\begin{pf}
We will prove Theorem~\ref{ETES} relying on the use of Theorem~\ref{LYAPU} for $\mathcal{A} = A+ BKC$, that generates, thanks to Assumption 1 and 2, an exponentially stable $C_0$-semigroup  on  $\mathcal{H} = H$ with
$$ \|e^{t(A+BKC)}z_0\|_H \leqslant M\|z_0\|_H e^{-\alpha t}\quad t\in \mathbb R_+.$$
For the same reason, hypothesis~\eqref{hypmax} is also satisfied with $c_1 \leqslant \|BKC \|_{\mathcal{L}(H)}$ and we can thus apply Theorem \ref{LYAPU}. Hence,  for $\beta \in (0, \alpha)$, there exists a Lyapunov functional $V_\beta : H \to \mathbb R_+^*$ satisfying $(i)$, $(ii)$ and $(iii)$ with constants $C_0$, $C_1$ and $C_2$ as in \eqref{Constants-eta}.
 


Let us now study the time derivative of $V_\beta$ along the trajectories (the state $z(t)$) of the event triggered system~\eqref{SystemEch}-\eqref{ETL1}.
Denoting $e_k(t)$ the measured deviation, we recall that the event-triggered law \eqref{ETL1} brings, for all $t \in [t_k,t_{k+1})$,
$$ \left\| e_k(t) \right\|_Y = \left\| C\left(z(t) - z(t_k)\right) \right\|_{Y} \leqslant \gamma \| z(t)\|.$$
Using successively statements $(ii)$, $(iii)$ and $(i)$ of Theorem \ref{LYAPU}, and the regularity properties proved in Proposition \ref{MoReg}, for any $z_0 \in \mathcal{D}(A)$ and $t\in [t_k,t_{k+1})$, we obtain
\begin{eqnarray*}
&\frac d{dt}&\left( V_\beta(z(t)) \right) =  \langle \nabla V_\beta(z(t)), \dot z(t)  \rangle_H 
\\
&=&  \langle \nabla V_\beta(z(t)), (A + BKC) z(t) \rangle_H  
\\
&&+ ~ \langle \nabla V_\beta(z(t)),  B K C( z(t_k) - z(t)) \rangle_H 
\\
&\leqslant& - 2 \beta V_\beta(z(t)) + C_0 \sqrt {V_\beta(z(t))} \| B K \|_{\mathcal L(Y,H)} \|e_k(t)\|_Y
\\
&\leqslant& - 2 \beta V_\beta (z(t)) + C_0 \gamma \sqrt {V_\beta(z(t))} \| B K \|_{\mathcal L(Y,H)}  \| z(t)\|_H\\
&\leqslant& - \left(2 \beta - \frac{C_0 \gamma}{\sqrt {C_1}}  \| B K \|_{\mathcal L(Y,H)}\right) V_\beta(z(t)).
\end{eqnarray*}
Thus, with $C_0$ and $C_1$ as in \eqref{Constants-eta}, defining 
$$
	\delta_\beta = \beta - \gamma \sqrt{2(1+ (c_1+ \beta)M^2)} \| BK \|_{\mathcal L(Y,H)}
$$
we get  
$$
	\frac d{dt} \left( V_\beta(z(t)) \right) 
	\leqslant - 2 \delta_\beta V_\beta(z(t)). 
$$
In particular, if $\delta_\beta$ is positive, we can write that for all $t \in (t_k, t_{k+1})$, 
$$
	 V_\beta(z(t)) \leqslant e^{-2\delta_\beta (t -t_k)} V_\beta (z(t_k)).
$$
Iterating on each interval $(t_k, t_{k+1})$, we easily get that for all $z_0 \in \mathcal{D}(A)$ and $t \geqslant 0$, 
$$
	 V_\beta(z(t)) \leqslant e^{-2\delta_\beta t} V_\beta (z_0), 
$$
which can be rewritten (recall $(i)$ and the values of $C_1$ and $C_2$ in \eqref{Constants-eta})
$$
	\| z(t) \|_H^2 \leqslant  \left( 1+ \frac{(c_1+\beta) M^2}{\alpha - \beta}\right) e^{-2\delta_\beta t} \| z_0\|_H^2.
$$
Since $c_1 \leq \| BK C\|_{\mathcal{L}(H)}$ and $\beta \in (0, \alpha)$, writing 
$$\tilde \delta_\beta =  \beta - \gamma \sqrt{2(1+ ( \| BK C\|_{\mathcal{L}(H)}+ \alpha)M^2)} \| BK \|_{\mathcal L(Y,H)}$$
we have $ \delta_\beta \geq \tilde \delta_\beta.$
Accordingly, if $\gamma$ satisfies \eqref{Cond-SmallGamma} and $\delta$ satisfies \eqref{Cond-delta}, there exists $\beta \in (0, \alpha)$ such that $\tilde \delta_\beta >\delta$, hence $\delta_\beta > \delta$, and the results follows.
$\hfill$ \qed 
\end{pf}

\begin{rem}
	When $A$ generates a unitary group on $H$, $B \in \mathcal{L}(U,H)$ and $A- BB^*$ generates an exponentially stable $C_0$ semigroup in $H$, all the above theory applies (corresponding to $K = - Id_U$, and $C = B^* \in \mathcal{L}(H,U)$) and the constant $c_1$ in \eqref{hypmax} can be taken as $c_1 = 0$. Accordingly, following the above proof of Theorem \ref{ETES}, the conditions on $\gamma$ and $\delta$ can be slightly relaxed. 
\end{rem}

\section{Examples of application}\label{Sec-Example}

The interest of our approach is that our abstract model encompasses many PDE models existing in the literature, such as transport equations with periodic boundary conditions, wave equation, Schr\"odinger equation, linear Korteweg-de-Vries equation, or linear beam models. We will give below some examples fitting our abstract framework and in which our results apply.

\paragraph*{A family of transport equations -}
Let us consider the following infinite dimensional system, consisting in a family of transport equations for speeds $v$ varying between $1$ and $\infty$:  Setting $\Omega = (0, 1) \times (1, \infty)$ and $\omega = (0,1/2)$, we consider the equation of solution $f=f(t,x,v)$ 
\begin{equation}
	\label{Eq-Family-Transport}
	\left\{
	\begin{array}{ll}
		\partial_t f + v \partial_x f + f {\bf 1}_{\omega}(x) = 0, \ &t \geqslant 0,\, (x,v) \in \Omega, 
		\\
		f(t,0,v) = f(t, 1, v), & t \geqslant 0, \, v \in (1, \infty), 
		\\
		f(0, x, v) = f_0 (x,v) & (x,v) \in \Omega.
	\end{array}
	\right.
\end{equation}
This system fits the abstract framework presented in this article with $H = L^2(\Omega)$, 
\begin{eqnarray*}
&&A f(x,v)= - v \partial_x f(x,v), ~ \hbox{ for } f\in  \mathcal{D}(A)~ \hbox{ where }\\
&&\mathcal{D}(A) = \big\{ f= f(x,v) \in H^1_x(0,1; L^2_v(1, \infty)) \\
&&\qquad\qquad\qquad\qquad\qquad\qquad \text { with } f(0, \cdot) = f(1, \cdot) \big\}, 
\end{eqnarray*}
$B = C = -K = {\bf 1}_{x \in \omega}$ and $U = Y = H$ for which it is then clear that Assumptions 1 and 2 are satisfied. To prove Assumption 3, let us remark that the system \eqref{Eq-Family-Transport} is a family of transport equations indexed by the velocity parameter $v$ on $(0,1)$ with periodic boundary conditions (i.e. on the torus $\R/\mathbb{Z}$). Therefore, for all $(t,x, v) \in (0, \infty) \times \Omega$, solutions of \eqref{Eq-Family-Transport} are given explicitly by the formula 
$$
	f(t, \{ x+ t v\} , v) = f_0(x,v) \exp\left( - \int_0^t \!1_\omega ( \{ x + s v\}) \, ds \right),  
$$
where $\{ \cdot \}$ here denotes the fractional part. Consequently, $t \mapsto \| f(t, \cdot, v)\|_{L^2(0,1)}$ decays and for all $x \in (0,1)$, $f(1/v, x, v) = f_0(x,v) \exp(- 1/(2v))$, so that we can deduce for all $t \geqslant 0$ and $v \geqslant 1$, using $\frac Nv \leq t \leq \frac{N+1}v$,
\begin{align*}
	\| f(t, \cdot, v) \|_{L^2(0,1)} 
	& \leqslant \left\| f\left(\frac Nv, \cdot, v\right) \right\|_{L^2(0,1)} 
	\\
	&\leqslant e^{- N/(2v)}\| f_0(\cdot, v)\|_{L^2(0,1)}\\
	& \leqslant  e^{ - t /2} e^{1/(2v )} \| f_0(\cdot, v)\|_{L^2(0,1)} 
	\\
	&\leqslant e^{(1-t)/2}\| f_0(\cdot, v)\|_{L^2(0,1)}.
\end{align*}
Taking the square and integrating in $v \in L^2(1, \infty)$, we immediately deduce Assumption 3 for \eqref{Eq-Family-Transport}, with $M = e^{1/2}$ and $\alpha = 1/2$.

Our results thus apply in this case, yielding that, provided $\gamma>0$ is small enough and $y_0 \in \mathcal{D}(A)$, the event-based law \eqref{discretecontrol}--\eqref{ETL}, which consists here in solving
$$
	\partial_t f + v \partial_x f + f(t_k, x,v) {\bf 1}_{\omega}(x) = 0, 
$$
for $ t \in [t_k, t_{k+1})$ and $(x,v) \in \Omega$ instead of the first equation of\eqref{Eq-Family-Transport}, generates an exponentially stable solution. 

Note that system \eqref{Eq-Family-Transport} is also a case in which no periodic sampling control stabilizes the system exponentially. To be more precise, take $T_1 >0$ and the sequence $T_k = k T_1$ for $k \in \N$, and consider the solution $f$ of 
\begin{equation}
	\label{Eq-Family-Transport-Periodic-Sampling}
	\left\{
	\begin{array}{ll}
		\partial_t f + v \partial_x f + f (T_k) {\bf 1}_{\omega}(x) = 0, \quad & 
		\\
		& \hspace{-3cm} \text{ for } t \in [T_k,T_{k+1}),\, (x,v) \in \Omega,  
		\\
		f(t,0,v) = f(t, 1, v), &\hspace{-1cm} t \geqslant 0, \, v \in (1, \infty), 
		\\
		f(0, x, v) = f_0 (x,v) &\hspace{-1cm}(x,v) \in \Omega.
	\end{array}
	\right.
\end{equation}
Indeed, if system \eqref{Eq-Family-Transport-Periodic-Sampling} was exponentially stable, there would exist $C_*$ and $\delta >0$  such that any solution would satisfy $\| f(t) \|_H \leqslant C_* \exp(- \delta t) \|f_0 \|_H$ for all $t \geqslant 0$ and any $f_0 \in H$. If so, let us choose $k_* \in \N$ such that $C_* \exp( - \delta T_{k_*}) \leqslant 1/2$. 
We should thus have that for all $f_0 \in H$, $\| f(T_{k_*}) \|_H \leq \| f_0 \|_H/2$. Now, take $v_1 = 1/ T_1$ if $T_1 \leqslant 1$ or a multiple of $1/T_1$ which is larger than~$1$ if $T_1 >1$. Then for all $t \in (0, \infty)$, the characteristic $t \mapsto \{ 3/4 + t v_1\}$ does not meet the set $\omega$ at times $(T_k )_{k \in \N}$. 
By continuity of the maps $(x,v) \mapsto \{x + T_k v \}$, one can find a neighborhood $\mathcal{N}$ of $(3/4, v_1)$ in $\Omega$ such that for all $(x_0, v) \in \mathcal{N}$, the characteristics $t \mapsto \{x_0 + t v\}$ do not meet the set $\omega$ at times $(T_k )_{k \in \{1, \cdots, k_*\}}$. 
Therefore, if we take $f_0$ supported in $\mathcal{N}$, the solution $f$ of \eqref{Eq-Family-Transport-Periodic-Sampling} is given by $f(t, \{ x + tv \}, v) = f_0(x,v)$ for all $t \in [0,T_{k_*}]$ and $(x,v) \in \Omega$ and satisfies for all $k \in \{1, \cdots, k_*\}$, $f (T_k) {\bf 1}_{\omega} = 0$. This of course contradicts the exponential stability of \eqref{Eq-Family-Transport-Periodic-Sampling}.

\paragraph*{Wave equation -}
Let $ \Omega $ be a smooth open  bounded domain of $\mathbb R^d$, and $p \in L^\infty(\Omega; R_+)$. We consider the wave equation damped in an open subset $\omega$, of solution $w=w(t,x)$:
\begin{equation}\label{eqW}
\hspace{-0.3cm}\left\{  
\begin{array}{ll}
w_{tt} - \Delta w  + p(x) w  = -k{\bf 1}_{\omega}(x) w_{t}, & \hbox{ in } \mathbb R_+\times \Omega,\\
w(t,x) = 0, \quad &\hbox{ in } \mathbb R_+\times \partial\Omega,\\
(w(0,\cdot), w_t(0,\cdot)) = (w_0, w_1),  & \hbox{ in } \Omega.
\end{array}
\right. 
\end{equation}
Here $k>0$ is the gain. We claim that this damped wave equation \eqref{eqW} fits into the abstract setting \eqref{System}, with 
\begin{align*}
&z = \begin{pmatrix} w\\ w_t\end{pmatrix}, \ 
A =  \begin{pmatrix} 0 & I \\ \Delta - p& 0 \end{pmatrix}, 
\\
&B=  \begin{pmatrix} 0 \\ {\bf 1_{\omega}}\end{pmatrix}, \ C = {}^t B,\ K = - k Id,
\end{align*}
$H = H^1_0(\Omega) \times L^2(\Omega) $, $\mathcal D(A) = H^2(\Omega) \cap H^1_0(\Omega) \times H^1_0(\Omega)$, and $U = Y = L^2(\Omega)$.

It is easy to check that $A^* = -A$ on the Hilbert space $H$ endowed with its canonical scalar product. Moreover, $B$, $C$ and $K$ are bounded linear operators. 

Hence, as soon as one can give a result proving the exponential stability of the closed loop system \eqref{eqW}, Theorem~\ref{ETES} proves subsequently the exponential stability of the corresponding event-triggered system \eqref{SystemEch}-\eqref{ETL1}. 

It is well-known that the exponential stability of the closed loop (damped) system \eqref{eqW} is equivalent to the observability of the (undamped) wave equation (see e.g. \cite{Haraux}), which is basically equivalent to the well-known Geometric Control Condition (see \cite{RauchTaylor1974,Bardos,BurqGerard}), stating that there exists $T>0$ such that all rays of geometric optics in $\Omega$ meet the domain $\omega$ in time less than $T$. Accordingly, under this same Geometric Control Condition, the event-triggered system \eqref{SystemEch}-\eqref{ETL1} corresponding to \eqref{eqW} is exponentially stable in the energy space $H = H^1_0(\Omega) \times L^2(\Omega) $ for initial data in $\mathcal D(A) = H^2(\Omega) \cap H^1_0(\Omega) \times H^1_0(\Omega)$. 
Such result obviously generalizes \cite{KBT-Auto22}, in which the damping was assumed to be acting on the whole set $\Omega$.

On a very similar setting, concerning now the Schr\"odinger equation set in a bounded domain, one can recover the result of \cite{KBT-SchroECC22} where an in-domain localized damping $-k {\bf 1}_{\omega}(x)w(t,x)$ is event-triggered. As for the wave equation, the Schr\"odinger equation with such a control fits into the abstract setting \eqref{System}, with $A = -i \Delta$ and $ B K C=  -k {\bf 1_{\omega}}$, when some appropriate geometric conditions over $\omega\subset \Omega$ ensuring stabilizability (equivalent to observability in this case) are satisfied (see for instance \cite{MachtZuazua,Lebeau1992} for classical sufficient geometric conditions, or \cite{Anantharaman-Leautaud,Burq-Zworski-2019} for more involved conditions). 

\paragraph*{Linear Korteveg-de Vries  equation -}
Let us consider the following linear Kortweg-de Vries (KdV) equation with periodic boundary conditions :
$$
\left\{  
\begin{array}{ll}
u_{t} (t,x) + u_{xxx}(t,x) = - k Gu(t,x),~~& \hbox{in } \mathbb R_+\times (0,2\pi  )\\
u(t,2\pi  ) = u(t,0), &  t\in \mathbb R_+\\
u_x(t,2\pi  ) = u_x(t,2\pi  ), &  t\in \mathbb R_+\\
u_{xx}(t,2\pi  ) = u_{xx}(t,0), &  t\in \mathbb R_+\\
u(0,x) = u_0(x), & x\in (0,2\pi  ) 
\end{array}
\right. 
$$
where the distributed control takes the form of the bounded operator defined by 
$$Gu(t,x) = g(x)\left(u(t,x) - \int_0^{2\pi}   g(s)u(t,s)\,ds \right)$$
for a given piecewise continuous non-negative function $g$ defined for $x\in [0,2\pi  ]$ with 
$$[g] = \int_0^{2\pi}   g(x)\,dx = 1 ~\hbox{ and  ~supp}(g) = [a,b] \subset [0,2\pi  ].$$
Such a control can be called a volume-conserving control given in feedback form, as presented in 
 \cite{RussellZhang-SICON93} where this system was studied first. It's interesting to see that it  corresponds to System \eqref{System} with 
\begin{align*}
&A =  - \partial_{xxx}, \quad
B K C=  -k G,
\quad H =  L^2(0,2\pi),
\\
&\text{ and } \mathcal D(A) = H^3_p(0,2\pi  ) = \big\{ u \in  H^3(0,2\pi  ), 
\\
& u(0) = u(2\pi ),  
u_x(0) = u_x(2\pi ) , 
u_{xx}(0) = u_{xx}(2\pi ) \big\}.
\end{align*}
It is easy to verify that $A$ generates a strongly continuous unitary group on $L^2(0, 2\pi)$ and that $-kG$ is a bounded control operator. Theorem 2 in \cite{RussellZhang-SICON93} states the exponential stabilization of the unique solution $u$ of KdV around the constant state equal to the total mass $[u_0] = \int_0^{2\pi} u_0(x) \, dx$ of the initial data $u_0$: $\forall t \geq 0$, 
$$
	\left\| u(t) - [u_0]\right\|_{L^2(0,2\pi)} \leqslant M e^{-\alpha t}\left\| u_0 - [u_0]\right\|_{L^2(0,2\pi)}.
$$
Although this is not directly in the framework of our result, we point out that, for solutions $u$ of KdV, the total mass $[u(t)] = \int_0^{2\pi} u(t,x) \, dx$ is independent of the time variable and thus is constant equal to $[u_0]$; and that $z(t) = u(t) - [u(t)] = u(t) - [u_0]$ solves the KdV equation (note that $Gu = Gz$) and from the above result, is exponentially stable. Therefore, to fit our abstract framework, we should consider the KdV equation on 
$$H_0 = L^2_0(0, 2\pi) =\left \{ z \in L^2(0,2\pi), \ \int_0^{2\pi} z(x) = 0\right\}$$ 
corresponding to the operator $A_0 = A$ with domain $\mathcal{D}(A_0) = \mathcal{D}(A) \cap H_0$, for which exponential stability holds. 

Theorem~\ref{ETES} then proves the non trivial exponential stability of the corresponding event-triggered system where the source term is now  $- k Gz(t_k,x)$ with appropriate law for triggering the $t_k$. Note that we did not make precise so far the choices of $B$, $K$ and $C$ such that $BKC = - k G$, and many choices are possibles, for instance $B = Id$, $K = -k Id$ and $C = G$. 
Consequently, we will recover the exponential convergence 
to the constant $[u_0]$ for solutions of the KdV equations with the event-triggered law given by 
\begin{multline*}
	t_{k+1} = \sup \Big\{ t>t_k, 
	\forall \tau \in [t_k,t), \\ 
	\left\| G\left(u(\tau) - u(t_k)\right) \right\|^2_{Y} \leqslant \gamma^2 \| u(\tau) - [u_0] \|_H^2\Big\}.
\end{multline*}

Let us finally point that a similar 1-d nonlinear KdV equation defined on a bounded domain was also studied for the same kind of event triggering control problematic in \cite{KBF}. There, the considered boundary conditions make the generator of the semigroup non self-adjoint and thus this case does not fit into our setting.  
Also, the work \cite{KBF} explores a scenario where the feedback is not only sampled over time, but also piecewise averaged in space, and regulated by a more technical event-triggering mechanism to maintain the stability of the system.

\section{Conclusion}

This article has proposed a unified setting that allows to encompass various closed loop systems of partial differential equations for which an event-triggering mechanism that samples the control input maintains the exponential stabilisation of the state.
The PDE has to be defined by a skew-adjoint operator and controlled and observed through bounded linear operators, and the continuously controlled closed loop system is assumed to be exponentially stable. 
The achievement of this paper was to prove that a simple and well-designed event-triggering law allows to keep this stability property, relying on the existence of the appropriate Lyapunov functional and ruling out any Zeno behavior. 

However, the proposed framework comes with several limitations. Most notably, it does not apply to parabolic equations such as the heat equation, due to the fact that these systems are not governed by skew-adjoint operators. Moreover, a key requirement of our approach is the boundedness of the operator $BKC$, ruling out boundary control setting for instance.

Another restriction lies in the linear nature of the considered systems: our results do not extend to nonlinear dynamics, which would require a more dedicate analysis. Additionally, to avoid Zeno behavior, a crucial assumption is that the continuous-time system must not decay faster than exponentially, a condition highlighted in Theorem~\ref{WellPo} \eqref{estimWP}, and central to the validity of our triggering mechanism.

Importantly, our sampling strategy does not rely on the explicit knowledge of a Lyapunov functional. Nevertheless, if such a functional were available, the same methodology could be applied, potentially leading to simpler or more flexible implementations. 

Future research directions include the development of dynamic event-triggering mechanisms, as recently proposed in the context of hyperbolic systems in \cite{WangKrsticTAC22} or for the Schrödinger equation in \cite{KBT-Dynamic22}. Another promising direction involves incorporating Luenberger observers into the event-triggered control design, in the spirit of the observer-based strategies developed in \cite{Espitia-SCL20} for boundary-controlled hyperbolic systems.

\bibliographystyle{alpha}
\bibliography{BiblioETL}

\end{document}